\documentclass[12pt]{article}
\usepackage{amssymb}
\usepackage{amsthm}\usepackage{epsf,epsfig,amsfonts}
\usepackage{amsmath}
\usepackage{graphicx}
\usepackage{a4wide}
\newtheorem{corollary*}{Corollary}

\newcommand{\be}{\begin{equation}}
\newcommand{\ee}{\end{equation}}

\newcommand{\R}{\mathbb{R}}

\newcommand{\weg}[1]{}

\usepackage{epsf,epsfig,amsfonts,a4wide}

\usepackage{amsfonts}
\usepackage{amsmath,amssymb,amsthm}
\usepackage{url}
\usepackage{amsmath,amssymb,amsthm}
\usepackage{amsfonts}
\usepackage{amsmath,amssymb,amsthm}
\usepackage{url}

\newtheorem{Th}{Theorem}

\newtheorem{Lemma}{Lemma}
\newtheorem{Cor}{Corollary}

\theoremstyle{remark}

\newtheorem{Rem}{Remark}
\newtheorem{Def}{Definition}

\newcommand{\eqdef}{\stackrel{{\rm def}}{=}}

\newcommand{\Z}{\mathbb Z}
\newcommand{\tr}{{\rm tr}}
\newcommand{\goth}{\mathfrak}

\newcommand{\Id}{\mbox{\rm Id}}

\newcommand{\const}{\mbox{\rm const}}
\begin{document}

\title{Fubini Theorem for pseudo-Riemannian metrics.}
\date{} \author{Alexey V. Bolsinov,
Volodymyr  Kiosak, 
Vladimir S. Matveev
}
\maketitle

 \begin{abstract}
 We generalize the following classical result of Fubini \cite{Fubini1}  for pseudo-Riemannian metrics: if three essentially  different  metrics on $M^{n\ge 3}$ share  the same unparametrized geodesics, and two of them (say, $g$ and $\bar g$) are strictly nonproportional (i.e., the minimal polynomial  of $g^{i\alpha } \bar g_{\alpha j}$ coincides with the characteristic polynomial) at least at  one point, then they have constant curvature.
\end{abstract}

\section{Introduction}
\subsection{Definitions and results}

Two Riemannian or pseudo-Riemannian metrics on the same manifold are said
to be \emph{geodesically equivalent}, if they have the same
geodesics considered as unparametrized curves.  Two metrics $g$ and
$\bar g$ are \emph{strictly nonproportional} at a point $x\in M$, if
the minimal polynomial  of the (1,1)-tensor $G^i_j:= g^{i k} \bar
g_{k j}$  coincides with the characteristic polynomial.  If one of
the metrics is Riemannian, strict nonproportionality  means that all
eigenvalues of $G$ have multiplicity 1.

The main result of the present paper  is the following

\begin{Th} \label{thm1}  Let $g$, $\bar g$ and $\hat g$ be three geodesically  equivalent metrics on a connected manifold $M^n $ of dimension $n\ge 3$. Suppose   there exists a  point  at which
$g$ and $\bar g$ are strictly nonproportional,  and    a point  at which
  $g$, $\bar   g$ and $\hat g$ are linearly independent. Then, the metrics $g, \bar g $ and $\hat g$  have constant curvature.
\end{Th}

If the metrics are Riemannian,  the local version of Theorem \ref{thm1} was proved by
  Fubini  in  \cite{Fubini1,Fubini2}.  The proof of Fubini is short and elegant, but unfortunately, does not work in the pseudo-Riemannian case. More precisely, Fubuni's proof is based on the following classical result:

\begin{Th}[Levi-Civita \cite{Levi-Civita}]
{\it  Assume $g$ and $\bar g$ are Riemannian geodesically equivalent
metrics that are strictly nonproportional at a point $p$.  Then, in
a neighborhood of this point there exists a coordinate system
$(x_1,...,x_n)$ such that the metrics $g,\  \bar g $ are  as
follows: \label{Levi-Civita}

\begin{eqnarray} \label{lcg}
ds^2_g &= & \sum_{i=1}^n  \left[\left|\prod_{\stackrel{j=1}{ j\ne i}}^n\left(X_i(x_i)- X_j(x_j)\right)\right| dx_i^2  \right], \\
ds^2_{\bar g} &= & \sum_{i=1}^n  \left[\frac{1}{X_i(x_i)\prod_{ \alpha = 1 }^n X_\alpha(x_\alpha)} \left|\prod_{\stackrel{j=1}{ j\ne i}}^n \left(X_i(x_i)- X_j(x_j)\right)\right| dx_i^2  \right],  \label{lcbg}
\end{eqnarray}
where  $X_i$ is a
  positive function of the
variable $x_i$ only. Moreover, every two metrics of this form are geodesically equivalent.     } \end{Th}

Fubini calculated the  riemannian curvature $R_{ij,km}$ of the metric \eqref{lcg} and observed that  the component $R_{ij,km}$ is non-zero only if $(i,j)= (k,m)$ or if $(i,j)= (m,k)$. Nowadays,  such curvature tensors are called \emph{diagonal}, see for example \cite{zaharov}.  Then he showed that, unless the curvature is constant,   the curvature tensor determines the coordinates lines $(x_1,...,x_n)$ uniquely. Thus, the metrics $g$, $\bar g$ and $\hat g$ must be simultaneously diagonalisable in a certain coordinate system and a short analysis shows that they are    linearly dependent at  every point.

Our proof, which is also valid in the pseudo-Riemannian setting, is
as follows. We study the partial differential equation \eqref{basic}
which is responsible for the fact that the metric $\bar g$ is
geodesically equivalent to $ g$. The unknown functions in this
equation are the components of a $(0,2)$-tensor $a_{ij}$ canonically
constructed by metric $\bar g$.   Then,   we find the integrability
conditions for this equation. This is a system of linear equations
on $a$ and the Hessian of $\lambda:= \tfrac{1}{2}a_{\alpha \beta }
g^{\alpha \beta}$; the coefficients  in these equations are algebraic
expressions of entries of the curvature tensor and the metric. We show that if the system has two solutions, then for certain $K\in \mathbb R$
the components of the curvature must satisfy the
condition $a_{i}^{ \alpha} Z_{\alpha jkl} +  a^{\alpha}_{j}Z_{\alpha ikl}=0$,
where $Z_{ijkl}= R_{ij,kl}- K \cdot ( g_{il}g_{jk} -g_{ik}g_{jl} )$. Then, it is an easy exercise in linear algebra to show that if in addition  the metrics $g$ and $\bar g$ are strictly nonproportional, then this  condition  on the curvature implies that the curvature is actually constant. Then, the curvature of the metrics $\bar g$ and $\hat g$ is  constant as well by the Beltrami Theorem.

\begin{Rem} We emphasize that the essential part of our proof is, in fact,  the analysis of algebraic properties of the integrability conditions  for  the  ''geodesic equivalence equation"  \eqref{basic}.  Algebraic aspects  of  the integrability conditions will  also be clarified in Section \ref{bolsinov} where we discuss an  unexpected relationship between geodesically equivalent metrics and the so-called sectional operators on semisimple Lie algebras.  As an application of this observation, we give  an alternative,  pure algebraic  proof of the (local version of the) Fubini theorem,  and we believe that  this idea  might be useful in a wider context.

\end{Rem}

\begin{Rem}  All assumptions in Theorem \ref{thm1} are important. Indeed:
\begin{itemize}
\item if the dimension $n$ is 2,  counterexamples were constructed in \cite{koenigs}, see also \cite{bryant};

\item if the metrics are not strictly nonproportional, counterexamples can be found in \S9 of  \cite{s1}, see also  \cite{shandra,s2,s3};

\item if we omit the
assumption that $g$, $\bar g$, and $\hat g$ are linearly
independent, then  we can take $\hat g= \bar g$, and  the
Levi-Civita metrics (\ref{lcg}) and (\ref{lcbg})  for generic $X_i$
give a counterexample to Theorem~\ref{thm1}.

\end{itemize}
\end{Rem}

\subsection{Motivation  }

The first  motivation, which was the reason why Fubini studied this question,  came from the study of projective vector fields of pseudo-Riemannian metrics. Recall that a vector field $v$ is \emph{projective} \weg{(\emph{affine}, respectively)}   with respect to the metric $g$ if its (local) flow takes geodesics to geodesics\weg{(preserves the Levi-Civita connection of $g$, respectively)}.   Projective and affine vector fields are  very classical objects  of investigation: both Lie \cite{Lie} and Schouten \cite{schouten}
explicitly formulated the problem of constructing all metrics admitting one or sufficiently many projective vector fields.
As a direct corollary of Theorem ~\ref{thm1} we obtain

\begin{Cor} \label{cor1} Suppose  that $g$ and  $\bar g$ are  geodesically equivalent metrics of a nonconstant curvature  on a connected manifold $M^n $ of dimension $n\ge 3$. Suppose   there exists a  point such that
$g$ and $\bar g$ are strictly nonproportional.    Then, the dimension of the space of projective vector fields minus the dimension of the space of homothety    vector fields is  at most one.
\end{Cor}

Recall that the vector field $v$ is \emph{ a homothety vector field,} if $L_v g = \const\cdot g$, where $L_v$ is the Lie derivative. We allow the case $\const =0$, so  Killing vector fields are also homothety vector fields.

\begin{Rem} In   dimension  2, Corollary \ref{cor1}  is wrong. Indeed, Darboux-superintegrable metrics  admit three projective vector fields and, as a rule,  only one homothety  vector field, which is a Killing vector field. The definition of Darboux-superintegrable metrics  and a description of their projective vector fields can be found in \cite{bryant}.   
Moreover, it is possible to show that there are no other counterexamples to the 2-dimensional version  of Corollary \ref{cor1}. This  is a very nontrivial statement which follows  from  the results of  \cite{bryant, alone}, where all 2-dimensional metrics admitting projective vector fields were constructed.
\end{Rem}

Moreover, it is possible to use Theorem ~\ref{thm1} in order to describe all projective vector fields of a metric $g$  under the assumptions that there exists a metric $\bar g$
which is geodesically equivalent to $g$ and strictly nonproportional to $g$.
Indeed, in this case the system of PDE on $v$ is a system of ODE's,
which is much easier to analyse than a system of PDE's, see \cite{Aminova1, Aminova2, bryant, Fubini1,  obata, sol,CMH, archive,    alone}, where this ODE-system was obtained and   completely solved in particular cases.

One more  motivation comes from the theory of superintegrable
systems. Recall that a metric is called \emph{superintegrable}, if
the number of independent integrals of special form is greater than
the dimension of the manifold. There are different possibilities for
the special form of integrals; de-facto the most standard special
form  of the integrals is  the so-called Benenti-integrals, which
are essentially the same as geodesically equivalent strictly
nonproportional metrics, see \cite{Benenti3,benenti,era1}.  The results
of the present paper show that    Benenti-superintegrable  metrics of
nonconstant curvature cannot exist, which was a folkloric
conjecture.

  \subsection{ History  }
The theory of  projective transformations has a long and fascinating
history. The first  non-trivial examples of projective transformations
were discovered  by   Lagrange \cite{lagrange,beltrami_short}.
Geodesically equivalent metrics were studied by Beltrami, Darboux,
Levi-Civita, Painleve,  Eisenhart, Weyl and other classics.  One can
find more historical details in the surveys
\cite{Aminova2,mikes,duna} and in the introduction to the papers
\cite{archive,hyperbolic,topology,dedicata,fomenko,sbornik,integrable,quantum}.

The pseudo-Riemannian version of Fubini's  Theorem  was investigated
by  the Kazan school of geometry, in particular  by A.\,Aminova and her
collaborators. In particular, they tried to mimic the Fubini proof for
pseudo-Riemannian metrics. They used the description of geodesically
equivalent pseudo-Riemannian metrics obtained by Aminova
\cite{Aminova}, which is a generalization of Levi-Civita's Theorem
\ref{Levi-Civita}. Unfortunately, the description by Aminova is very
complicated, which makes this programm to be computationally very
hard. Besides, there are infinitely many different types of normal
forms for pseudo-Riemannian metrics (depending on the Jordan form of
the tensor $G$), and each of these types require a separate
inverstigation. According to \cite{zakirova},  they proved
Theorem~\ref{thm1} in  dimensions up to 6, though we did not find
the place where  the proof is written.

\subsubsection*{Acknowledgements}
Alexey Bolsinov thanks Russian Foundation  for Fundamental Research  for partial financial support 
and Mathematisches Institut, Friedrich-Schiller-Universit\"at, Jena  for hospitality.
Vladimir Matveev   is grateful to   R.~Bryant and   M.~Eastwood for useful discussions, and MSRI for the hospitality.  Volodymyr Kiosak and Vladimir Matveev thank   
Deutsche Forschungsgemeinschaft (Priority Program 1154 --- Global Differential Geometry)
  and FSU Jena for partial financial support.

 \section{Proof}
 \subsection{ Schema of the proof} \label{schema}

We first show (Section \ref{tmp}) that if geodesically equivalent
metrics are strictly non-proportional at least at one point, than it
is so at almost every point.  Then, we will show (Section
\ref{tmp1}) that if $g$, $\bar g$, and $\hat g$  are linearly
dependent at every point of a neighborhood, they are linearly
dependent at every point of the manifold.

Combining these two observation, we see that if the assumptions of
Theorem ~\ref{thm1} are fulfilled at one point, they are fulfilled
at almost every point, i.e., though the theorem is global, it is
sufficient to prove it locally. This will be done in Sections
\ref{main}, \ref{111}.

 \subsection{What we will use in the proof}  \label{mittel}

There are many  tensor reformulations of the condition ``metrics $g$
and $\bar g $ are geodesically equivalent". In our paper we will use
the following one, which was suggested by Sinjukov \cite{Sinjukov},
see also \cite{benenti,eastwood}: given two metrics $g$ and $\bar g$,
consider the $(0,2)-$tensor
 \begin{equation} \label{a}
 a_{ij}:=  \left|\frac{\det(\bar g)}{\det(g)}\right|^{\frac{1}{n+1}}\cdot g_{i\alpha} \bar g^{\alpha \beta }  g_{j\beta}, \end{equation}
 and the function
 \begin{equation}\label{lambda}  \lambda := \frac{1}{2} a_{\alpha \beta} g^{\alpha \beta} \end{equation}
 where and $g^{\alpha\beta}$ and $\bar g^{\alpha\beta}$ denote the tensors dual to $g_{ij}$ and $\bar g_{ij}$ respectively, i.e.,
$g^{i\alpha} g_{\alpha j}= \delta_j^i$ and $\bar g^{i\alpha} \bar
g_{\alpha j}= \delta_j^i$.
 \begin{Th}\cite{Sinjukov, benenti}
 The metrics $g$ and $\bar g$ are geodesically equivalent, if and only if
 \begin{equation} \label{basic}
 a_{ij,k}= \lambda_{,i} g_{jk} + \lambda_{,j}  g_{ik},
 \end{equation}
 where the covariant derivative is taken with respect to the
 (Levi-Civita connection of the) metric $g$.
 \end{Th}

 We will also use the following connection between geodesically equivalent metrics and integrable geodesic flows due to \cite{MT}.

Let $a_j^i$ be as in \eqref{a}, with one index lifted by $g$ (so now $a$ is a $(1,1)-$tensor,
   self-adjoint  with respect to $g$).
Consider the family $S_t$ of $(1,1)$-tensors
\begin{equation}\label{st}
 S_t\eqdef \det(a - t\ \mbox{\rm Id})\left(a-t\ \mbox{\rm Id}\right)^{-1}, \ \  \ t\in \mathbb{R}.
 \end{equation}
\begin{Rem}
Although $\left(a-t\ \Id\right)^{-1}$ is not defined for
$t$
lying
in the spectrum of $a$, the tensor  $S_t$  is well-defined
for every  $t$. Actually, the matrix of $S_t$ is the comatrix of the matrix of $a- t \ \Id$. In particular,   $S_t$
is a  polynomial
 in $t$ of degree $n-1$
with coefficients being  (1,1)-tensors.
\end{Rem}
\noindent We will identify the tangent and  cotangent bundles of $M^n$ by $g$.
This identification allows us to transfer the
 natural  Poisson structure from  $T^*M^n$ to   $TM^n$.

\begin{Th}[\cite{ MT, Topalov, dedicata}]
\label{integrability}
 If $a$ is a solution of \eqref{basic},  then  for every  $t\in \mathbb{R}$  the function
\begin{equation}\label{integral}
I_{t}:TM^n\to \mathbb{R}, \ \ I_{t}(\xi)\eqdef g(S_{t}(\xi),\xi)
\end{equation}
is an  integral for the geodesic flow
 of  $g$.
\end{Th}

Recall that a  function is an \emph{integral} of the geodesic flow, if it is constant  along the  trajectories, i.e., along the curves on the tangent bundle of the form $(\gamma(t), \dot \gamma(t))$, where $\gamma$ is a geodesic and $\dot\gamma(t)$ is its velocity vector at a point $t$.

We will also  use the following statement, whose Riemannian version was obtained in  \cite{MT}, and pseudo-Riemannian generalisation is  due to Topalov \cite{involutivity}. As in the introduction, we denote by $G$ the $(1,1)-$tensor $G^i_j:= g^{i k} \bar
g_{k j}$.  

\begin{Th} [follows from  Theorem 2 and Section III B  of \cite{involutivity}] \label{4}
Suppose the degree of the minimal polynomial of $G$  is $r$ at every point of the neighbourhood $U(p)$. Then, for arbitrary mutually  different  $t_1,...,t_{r+1}\in \mathbb{R}$ the functions $I_{t_1}$,...,$I_{t_r}$ are functionally independent, and
the function $I_{t_{r+1}}$ is a linear combination of the
functions $I_{t_1}$,...,$I_{t_r}$ (the coefficients of the linear combination are constant).
\end{Th}
Recall that functions $f_1,...,f_r$ are \emph{functionally independent}, if their differentials are linearly independent almost everywhere.

In other words, the number of independent integrals among  $I_t$ is
the degree of the minimal polynomial. In particular, if the metrics
are strictly nonproportional at a point $p$ (which immediately
implies that  they are strictly non-proportional at every point
in a small neighbourhood $U(p)$), then the differentials of the
integrals are linearly independent at almost every point of $TU(p)$.

\subsection{If geodesically equivalent metrics are strictly nonproportional at one point, then they are strictly nonproportional at almost every point. }  \label{tmp}

Let geodesically equivalent metrics $g$ and $\bar g$ be  strictly
nonproportional at $p$. Consider a  geodesic $\gamma$ passing
through $p:= \gamma(0)$. Let us show that  every point $q:=
\gamma(\tau)$ of the geodesic   has a neighborhood $U(q)$ such
that at almost every point of the neighborhood the metrics are
strictly nonproportional.  Since every point can be reached by a
finite sequence of geodesics, and since the condition ``the minimal
polynomial of $G$ has degree $n$" is an open condition,  this will
prove the statement formulated in the title of this section.

As we recalled in Section \ref{mittel}, for almost every point of
$TU(p)$ the differentials of the integrals $I_{t_1},..., I_{t_n}$
are linearly independent. Take a sequence of points $(p_k, v_k)\in
TU(p)$ converging to $(\gamma(0), \dot\gamma(0))$ such that 
$dI_{t_1}(p_k, v_k),\dots , dI_{t_n}(p_k, v_k)$ are linearly independent. Consider the
sequence of the geodesics $\gamma_k$ such that $\gamma_k(0)= p_k$
and $\dot\gamma_k(0)= v_k$. For sufficiently large $k$, the geodesics
exist up to time $\tau$ and the sequence
$(\gamma_k(\tau),\dot\gamma_k(\tau))$ converges to $(\gamma(\tau),
\dot\gamma(\tau))$.

Since the  integrals are preserved by the geodesic flow, the
differentials of the integrals are also preserved by the geodesic
flow. Then, at the points $(\gamma_k(\tau),\dot\gamma_k(\tau))$ the
differentials of $I_{t_1},..., I_{t_n}$  are linearly
independent. By Theorem \ref{4}, in an arbitrary small neighborhood
of $\gamma_k(\tau)$ there exists a point $q_k$ such that  the
metrics are strictly nonproportional at  $q_k$ implying the claim, \qed

\subsection{ If  $g_{i k}$, $ \bar  g_{k j}$, and $\hat  g_{k j}$
 are linearly dependent at every point of a neighborhood $U$, then they are linearly dependent
 at every point  of $M$.} \label{tmp1}

Within this section we assume that $g$, $\bar g$, $\hat g$ are geodesically equivalent metrics on a connected manifold $M$ of dimension $n\ge 3$.
We consider the tensor $a_{ij}$ and the function $\lambda$ given by (\ref{a},\ref{lambda}).
The same objects for the pair of metrics $g, \hat g$ will be denoted by the  capital letters $A$ and $\Lambda$, i.e.,
\begin{equation} \label{A}
 A_{ij}:=  \left|\frac{\det(\hat g)}{\det(g)}\right|^{\frac{1}{n+1}}\cdot g_{i\alpha} \hat g^{\alpha \beta }  g_{j\beta} , \ \ \  \  \  \Lambda := \frac{1}{2} A_{\alpha \beta} g^{\alpha \beta} \end{equation}

We will first prove the following statement (essentially due to Weyl \cite{Weyl}).

\begin{Lemma}\label{weyl} Suppose  $a$ and $A$ are solutions of \eqref{basic}. Assume
 $a = C \cdot A$, where $C$ is a function.  Then, $C$ is a constant.
\end{Lemma}

{\bf Proof. } Our proof is different from the proof of Weyl and is
based on the ideas developed in  \cite{dedicata}.  Note that in the
proof we use only the fact that the dimension is greater than one,
i.e., it works in  dimension 2 as well.

Consider the integrals $I_0$ from \eqref{integral} constructed by
$a$ (we keep the notation  $I_0$ for it) and by $\hat g$ (we denote it by
$\mathcal{I}_0$). If $a = C \cdot A$, then
 the integrals  $I_0$ and $\mathcal{I}_0$ are proportional as well, direct calculations show that
   $I_0(\xi)= \pm C^{n-1}\cdot  \mathcal{ I}_0(\xi)$. Since the  functions $I_0$ and $\mathcal{I}_0$ are constant along  every trajectory of the geodesic flow, the coefficient of proportionality of these functions  is also constant along every  trajectory of the geodesic flow implying  that  it is constant everywhere.  \qed

Now let us   assume that at every point of the neighborhood $U$  the
tensors  $g$, $\bar g ,$ and  $\hat g$ are linearly dependent. Then,
for certain functions $c , d $  the tensors  $g, a, A$ satisfy
(probably in a smaller neighbourhood $U'\subseteq U$; without loss
of generality we can think that $U' $ coincides with $U$.)
\begin{equation} \label{3}
a_{ij} = c \, g_{ij} + d \, A_{ij}.
\end{equation}

We will   show that  the functions  $c , d $  are actually
constants.

Differentiating \eqref{3}  and substituting  \eqref{basic} and its analog for the solution $A$, we obtain

\begin{equation} \label{a1}
\lambda_{,i} g_{jk} + \lambda_{,j} g_{ik} = c_{,k}\, g_{ij} + d \,
\Lambda_{,i} g_{jk} + d \, \Lambda_{,j}g_{ik} + d _{,k}\,
A_{ij},\end{equation} which is evidently equivalent to
\begin{equation} \label{a2}
\tau_{i} g_{jk} + \tau_j g_{ik} = c _{,k}g_{ij} + d
_{,k}A_{ij},\end{equation} 
where $\tau_i= \lambda_{,i} - d\,
\Lambda_{,i}$. We see that for every fixed $k$ the left-hand side of (\ref{a2}) is
a symmetric matrix  of the form $\tau_i v_j + \tau_j v_i.  $   If $c
_{,k}$ is not proportional to $d _{,k}$, this will imply that
$g_{ij}$ also is of the form $\tau_i v_j + \tau_j v_i,   $ which
contradicts the non-degeneracy of $g$.  If $c _{,k}=f\cdot d _{,k}$,
then the coefficient $f$  of the proportionality should be a constant
implying $d = \const\cdot  {c }$  and $a_{ij} = c \, (g_{ij} +
\const\cdot A_{ij})$. Since the equation \eqref{basic} is linear and
since $g_{ij} $ and $A_{ij}$ are its solutions, their sum $g_{ij} +
\const\cdot  A_{ij}$ is also solution of \eqref{basic}. Then,  Lemma
\ref{weyl}   implies  that $c $ is constant. Thus, $c $ and $d $ are
constant in a neighbourhood $U$. Since the equation \eqref{basic}
is linear and of finite type, see \cite{eastwood}, linear
dependence of solutions in a neighbourhood implies linear dependence
of the solutions everywhere.   \qed

\begin{Rem} Though we used that the dimension of the manifold is at least three, the statement is true in dimension two as well provided the curvature of $g$ is not constant, see \cite{kruglikov}.

\end{Rem}

 \subsection{ Main step of the proof of Theorem \ref{thm1}} \label{main}
As we explained in Section \ref{schema}, in view of Sections \ref{tmp}, \ref{tmp1}, it is sufficient to prove Theorem locally only.

Assume the metrics $g, \bar g, $ and $\hat g$ are linearly independent at a point $p$.
Assume  the  metrics $g, \bar g $ are strictly nonproportional at $p$.
Then, for a certain neighborhood $U(p)$
 the metrics $g, \bar g, $ and $\hat g$ are linearly independent, and the metrics $g, \bar g $ are strictly nonproportional  at every  point of $U(p)$.

 Then, the solution $a$ corresponding to $\bar g$,   the solution $A$ corresponding to $\hat g$, and the metric $g$ are linearly independent at every point of $U(p)$, and the minimal polynomial of the
 (1,1)-tensor $a_j^i$ has degree $n$.

As in Section \ref{tmp}  we denote by $\Lambda $ the function \eqref{lambda} constructed by $A$.

 Let us first consider the case when $\lambda$ is constant. In this case, the equation \eqref{basic} implies that $a_{ij}$ is covariant constant. Since $a^i_{j}$ is self-adjoint with respect to $g$,
 at every point
  there exists a  basis $(b_1(p), ..., b_n(p))$
   such that the matrix of $g$ is diagonal with $\pm 1$ on the diagonal
   and the matrix of $a^i_{j}$ is Jordan-matrix.
   The basis is unique up to signs of the vectors, so locally $b_1$ can be taken to be smooth vector fields on the manifold.   Since the vectors of the basis are invariantly constructed by two  covariantly-constant object, they  are covariantly constant as well.
   Then, the metric $g$ is flat implying that $\bar g$ and $\hat g$ have constant curvature.

 In what follows we assume that $\lambda$ is not constant.

 Theorem \ref{thm1} is an easy consequence of the following two lemmas.
 
\begin{Lemma} 
\label{13}   
Assume that $a$ and $A$ are solutions of \eqref{basic} such that $a,A,g$ are linearly independent. Then,
at  every point  there exists $K$ such that  the tensor  $Z_{ijkl}:= R_{ij,kl} -
K\cdot  \left( g_{il}g_{jk} -  g_{ik}g_{jl} \right)$ satisfy the condition
\begin{equation}  a_{i}^{ \alpha} Z_{ \alpha jkl} +  a^{\alpha}_jZ_{\alpha i kl}=0,  \label{int33b}
\end{equation}
\end{Lemma}

We see that by construction the tensor $Z_{ijkl}$ is skew-symmetric with respect to the first two indexes. Note that  in Lemma \ref{13} we did not assume that the minimal polynomial of $a$ has degree $n$, i.e., the lemma is valid in a slightly more general setting than we need it.

\begin{Rem}  Tensor $Z$ plays important role in the theory of geodesically equvialent metrics; it appears quite naturally in the investigation of special metrics such that Einstein metrics, pseudo-symmetric metrics, K\"ahler metrics,     see for example \cite{kiosak,kiosak-mikes} for details. 
\end{Rem}

The next lemma shows that if in addition   the minimal polynomial of $a$ has degree $n$, then the condition \eqref{int13} implies $Z_{ijkl}= 0$.
 Then, the curvature of $g$ is constant. By Beltrami's Theorem  (see for example \cite{beltrami_short}, or the original papers \cite{Beltrami} and \cite{schur}),  the metrics $\bar g $ and $\hat g$ have constant curvature as well which proves Theorem ~\ref{thm1}.

Moreover, we will not use the indexes $k$ and $l$ in the proof, so Lemma \ref{last} is the matrix reformulation of the condition  $a_{i}^{ \alpha} Z_{ \alpha j} +  a^{\alpha}_jZ_{\alpha i }$ for arbitrary $(1,1)-$tensor $a$ and  a skew-symmetric (0,2)-tensor $Z$.

\begin{Lemma} \label{last}
Let $Z$, $a$  be $n\times n$-matrices such that $Z$ is skew-symmetric. Assume
 the minimal polynomial of $a$ has degree $n$. If $Za - a^t Z=0$, then $Z=0$.
\end{Lemma}

{The proof of Lemma \ref{last} }  is an easy exercise in linear algebra and  will be  left to
 the reader (it
 is  a kind of problem such that it is easier to prove than to understand the proof.)
 We recommend to consider the coordinates such that the matrix
 $a$ is in   Jordan form, and then to  calculate the matrix  $Za$. One immediately sees that if the matrix $Za$ is symmetric, then $Z=0$. \qed

 \subsection{Proof of Lemma \ref{13}.}  \label{111}

The proof is straightforward tensor calculations. The geometry behind the calculation could be understood with the help of \cite{eastwood}. There, the equations \eqref{basic} were written in the projectively invariant form (so that the equations for $g$ and for the geodesically equivalent metric $\bar g$ are the same). The prolongation of the equations was written as a connection on the projective tractor  bundle, and the curvature of the connection  was  calculated.  Its first part  related to (the analog of the  tensor) $a$ is a trace-free  object.  Moreover, for most objects, the trace-free  part of the covariant derivative coincides with the trace-free part of the corresponding derivative on the
tractor bundle.  

We will see that in the proof we consider  the integrability conditions for the equations
\eqref{basic}, and then    "artificially" write all objects   in trace-free form.  At the end we obtain the required equation \eqref{int33b}.

Note that in the proof we  will essentially use the symmetries of the Riemannian curvature tensor, which have sense
 only if the affine connection is the Levi-Civita connection of a metric. That means, it is important for us that among the solutions of the projective-invariant analog of $a$ there is a non-degenerate solution.

In  Section \ref{bolsinov},  we give another proof of Lemma \ref{13} which is, in fact, pure algebraic and  based on some ideas from the theory of integrable Hamiltonian systems on Lie algebras.

  {\bf Proof of Lemma  \ref{13}.}
 Integrability conditions for the equation \eqref{basic} are
 (we use the standard fact that  $a_{ij,kl}- a_{ij,lk}= a_{i \alpha }R^{\alpha}_{jkl} +  a_{\alpha j}R^{\alpha}_{ikl}$  for any  $(0,2)-$tensor $a_{ij}$)

\begin{equation}  a_{i \alpha }R^{\alpha}_{jkl} +  a_{\alpha j}R^{\alpha}_{ikl} =\lambda_{, li} g_{jk}+\lambda_{, lj} g_{ik}-\lambda_{, ki} g_{jl}-\lambda_{, kj} g_{il}. \label{int1} \end{equation}

 The same is true for the other solution $A$
 \begin{equation}  A_{i \alpha }R^{\alpha}_{jkl} +  A_{\alpha j}R^{\alpha}_{ikl} =\Lambda_{, li} g_{jk}+\Lambda_{, lj} g_{ik}-\Lambda_{, ki} g_{jl}-\Lambda_{, kj} g_{il}. \label{int2} \end{equation}

Starting from this point, the proof is purely algebraic:  the statement of  Lemma  \ref{13} is an algebraic corollary of  (\ref{int1}) and (\ref{int2}). We emphasize this once again in the next section by giving another version of the proof in the Lie-algebraic language. 

 Now let us  multiply the equation \eqref{int1} by $A^{l}_{s}$ and sum over $l$.
After renaming indexes,   we obtain
\begin{equation}  a_{i \alpha }R^{\alpha}_{jk\beta} A^{\beta}_l +
a_{\alpha j}R^{\alpha}_{ik\beta} A^{\beta}_l =\lambda_{, \alpha i} A^{\alpha}_l   g_{jk}+
\lambda_{, \alpha j}  A^{\alpha}_l  g_{ik}-\lambda_{, ki}   A_{jl} -\lambda_{, kj} A_{il}. \label{int3} \end{equation}

 Using the symmetry of the Riemann tensor we obtain   $a_{\alpha j}R^{\alpha}_{ik\beta} = a^{\alpha }_i R_{\alpha j, k \beta} A^{\beta}_l = a_i^\alpha R_{\beta k ,j \alpha } A^{\beta}_l  =a_i^\alpha A_{\beta l} R^{\beta}_{k j \alpha}$. Substituting this in  \eqref{int3}, we get

\begin{equation}  a_{i}^\alpha A_{\beta l} R^{\beta}_{ k i\alpha}  +  a_{j}^\alpha A_{\beta l} R^{\beta}_{ k j\alpha}  = \lambda_{, \alpha i} A^{\alpha}_l   g_{jk}+
\lambda_{, \alpha j}  A^{\alpha}_l  g_{ik}-\lambda_{, ki}   A_{jl} -\lambda_{, kj} A_{il}. \label{int4} \end{equation}

Let us now symmetrise  \eqref{int4} by $l,k$ to obtain
\begin{equation}\begin{array}{l}
a_{i}^\alpha \left(A_{\beta l} R^{\beta}_{ k j\alpha}  +   A_{\beta k} R^{\beta}_{ l j\alpha} \right)+
a_{j}^\alpha\left( A_{\beta k} R^{\beta}_{ l i\alpha}  +   A_{\beta l} R^{\beta}_{ k i\alpha}\right)  \\
= \lambda_{, \alpha i} A^{\alpha}_l   g_{jk}+
\lambda_{, \alpha j}  A^{\alpha}_l  g_{ik}-\lambda_{, ki}   A_{jl} -\lambda_{, kj} A_{il}+
\lambda_{, \alpha i} A^{\alpha}_k   g_{jl}+
\lambda_{, \alpha j}  A^{\alpha}_k  g_{il}-\lambda_{, li}   A_{jk} -\lambda_{, lj} A_{ik}.\end{array} \label{int5}\end{equation}

We see that the components staying in brackets are the left-hand side of the equation
\eqref{int2} with other indexes. Substituting  \eqref{int2} in the first term of the left-hand side we obtain

\begin{equation}\begin{array}{cc}
 a_{i}^\alpha \left(A_{\beta l} R^{\beta}_{ k j\alpha}  +   A_{\beta k} R^{\beta}_{ l j\alpha} \right) & =
a_{i}^\alpha \left(\Lambda_{,\alpha l} g_{jk}  +   \Lambda_{, \alpha k} g_{jl } -
\Lambda_{,j  l} g_{k\alpha }  -   \Lambda_{, j  k} g_{l \alpha}  \right) \\&  = a^\alpha_i \Lambda_{,\alpha l} g_{jk} + a_i^\alpha \Lambda_{,\alpha k}g_{jl}- \Lambda_{,jl} a_{ik} -  \Lambda_{,jk} a_{il}.\end{array}
\label{int6}
\end{equation}
Similarly, the second term of the left-hand side is
\begin{equation}\begin{array}{cc}
 a_{j}^\alpha \left(A_{\beta l} R^{\beta}_{ k i\alpha}  +   A_{\beta k} R^{\beta}_{ l i\alpha} \right) & =
a_{j}^\alpha \left(\Lambda_{,\alpha l} g_{ik}  +   \Lambda_{, \alpha k} g_{il } -
\Lambda_{,i  l} g_{k\alpha }  -   \Lambda_{, i  k} g_{l \alpha}  \right) \\&  = a^\alpha_j \Lambda_{,\alpha l} g_{uk} + a_j^\alpha \Lambda_{,\alpha k}g_{il}- \Lambda_{,il} a_{jk} -  \Lambda_{,ik} a_{jl}.\end{array}
\label{int7}
\end{equation}
Substituting (\ref{int6},\ref{int7}) in \eqref{int5}, we obtain
\begin{equation}\begin{array}{l}
a^\alpha_i \Lambda_{,\alpha l} g_{jk} + a_i^\alpha \Lambda_{,\alpha k}g_{jl}- \Lambda_{,jl} a_{ik} -  \Lambda_{,jk} a_{il}+
a^\alpha_j \Lambda_{,\alpha l} g_{ik} + a_j^\alpha \Lambda_{,\alpha k}g_{il}- \Lambda_{,il} a_{jk} -  \Lambda_{,ik} a_{jl}  \\
= \lambda_{, \alpha i} A^{\alpha}_l   g_{jk}+
\lambda_{, \alpha j}  A^{\alpha}_l  g_{ik}-\lambda_{, ki}   A_{jl} -\lambda_{, kj} A_{il}+
\lambda_{, \alpha i} A^{\alpha}_k   g_{jl}+
\lambda_{, \alpha j}  A^{\alpha}_k  g_{il}-\lambda_{, li}   A_{jk} -\lambda_{, lj} A_{ik}.\end{array} \label{int8}\end{equation}

 Collecting the terms by $g$, we  see that \eqref{int8}  can be written as

 \begin{equation}\begin{array}{l}
\left(a^\alpha_i \Lambda_{,\alpha l} -  \lambda_{, \alpha i} A^{\alpha}_l  \right) g_{jk}
  + \left(a_i^\alpha \Lambda_{,\alpha k}-
\lambda_{, \alpha i} A^{\alpha}_k  \right) g_{jl}
+  \left(a^\alpha_j \Lambda_{,\alpha l}  -\lambda_{, \alpha j}  A^{\alpha}_l \right) g_{ik}
+ \left(a_j^\alpha \Lambda_{,\alpha k}-
\lambda_{, \alpha j}  A^{\alpha}_k \right) g_{il} \\
=
 \Lambda_{,jl} a_{ik} +  \Lambda_{,jk} a_{il} + \Lambda_{,il} a_{jk} +  \Lambda_{,ik} a_{jl}   -\lambda_{, ki}   A_{jl} -\lambda_{, kj} A_{il}-\lambda_{, li}   A_{jk} -\lambda_{, lj} A_{ik}.\end{array} \label{int9}\end{equation}

We denote $\tau_{il}:= a^{\alpha}_i \Lambda_{, \alpha l} - A^{\alpha}_l \lambda_{, \alpha i}$. In this notation,  the equation \eqref{int9}  is
\begin{equation}\begin{array}{lll}
&&\tau_{il}   g_{jk}
  + \tau_{ik}  g_{jl}
+  \tau_{jl} g_{ik}
+ \tau_{jk} g_{il} \\
&=&
 \Lambda_{,jl} a_{ik} +  \Lambda_{,jk} a_{il} + \Lambda_{,il} a_{jk} +  \Lambda_{,ik} a_{jl}   -\lambda_{, ki}   A_{jl} -\lambda_{, kj} A_{il}-\lambda_{, li}   A_{jk} -\lambda_{, lj} A_{ik}.\end{array} \label{int10}\end{equation}
 Let us show that $\tau$ is symmetric.
 Multiplying with $g^{jk}$ and contracting with respect to $j,k$,  we obtain
   \begin{equation}\begin{array}{lcl}
(n+2) \tau_{il}
+ \left(\tau_{jk} g^{jk}\right)  g_{il} &
=&
 \Lambda_{,\alpha l} a_{i}^\alpha  +  \left(g^{jk}\Lambda_{,jk}\right) a_{il} +
 \Lambda_{,il} \left(a_{jk}g^{jk} \right) +  \Lambda_{,i\alpha} a^{\alpha}_l   \\ && -\lambda_{, \alpha i}   A^{\alpha}_l -\left(\lambda_{, kj} g^{kj}\right) A_{il}-\lambda_{, li}   \left(A_{jk} g^{jk}\right) -\lambda_{, l\alpha } A_{i}^\alpha .\end{array} \label{int11}\end{equation}
We see that the right-hand side is symmetric w.r.t. $i,l$. Then, so should be the left-hand-side implying
$\tau_{il}=\tau_{li}$.
We also see that the sum of the  first, fourth,  fifth and last terms of the right-hand side is $\tau_{il}+ \tau_{li} = 2 \tau_{il}$. Then,  the equation \eqref{int11} is equivalent to
\begin{equation}
 \tau_{il}
=\frac{1}{n}\left(-\left(\tau_{jk} g^{jk}\right)   g_{il} +
 \left(g^{jk}\Lambda_{,jk}\right) a_{il} +
 \Lambda_{,il} \left(a_{jk}g^{jk} \right)
  -\left(\lambda_{, kj} g^{kj}\right) A_{il}-\lambda_{, li}   \left(A_{jk} g^{jk}\right)\right) . \label{int12}\end{equation}

 Now we return to the equation \eqref{int10}. We alternate the equation with respect to $j, k$.
 \begin{equation}\begin{array}{lll}
&&   \tau_{ik}  g_{jl}
+  \tau_{jl} g_{ik} -   \tau_{ij}  g_{kl}
-   \tau_{kl} g_{ij}
  \\
&=&
 \Lambda_{,jl} a_{ik}   +  \Lambda_{,ik} a_{jl}   -\lambda_{, ki}   A_{jl}  -\lambda_{, lj} A_{ik} - \Lambda_{,kl} a_{ij}   -  \Lambda_{,ij} a_{kl}   +\lambda_{, ji}   A_{kl}  +\lambda_{, lk} A_{ij}     . \end{array} \label{int13}\end{equation}

Let us now rename $i\leftrightarrow k$ in \eqref{int13}. We obtain
\begin{equation}\begin{array}{lll}
 &&  \tau_{ik}  g_{jl}
+  \tau_{jl} g_{ik} -   \tau_{kj}  g_{il}
-   \tau_{il} g_{kj}
  \\
&=&
 \Lambda_{,jl} a_{ik}   +  \Lambda_{,ik} a_{jl}   -\lambda_{, ki}   A_{jl}  -\lambda_{, lj} A_{ik} - \Lambda_{,il} a_{kj}   -  \Lambda_{,kj} a_{il}   +\lambda_{, jk}   A_{il}  +\lambda_{, li} A_{kj}     . \end{array} \label{int14}\end{equation}

 Adding  \eqref{int10} with \eqref{int14}  and dividing by $2$ for cosmetic
 reasons, we obtain

 $$
\tau_{ik}  g_{jl}
+   \tau_{jl} g_{ik}
=
 \Lambda_{,jl} a_{ik}  +  \Lambda_{,ik} a_{jl}
  -\lambda_{, ki}   A_{jl}      -\lambda_{, lj} A_{ik}    . $$

  Substituting the expression for $\tau$
  from  \eqref{int12}, we obtain

  $$\begin{array}{ccl}
&&\frac{1}{n}\left(-\left(\tau_{\beta\alpha } g^{\beta\alpha }\right)   g_{ik} +
 \left(g^{\alpha \beta }\Lambda_{,\alpha \beta }\right) a_{ik} +
 \Lambda_{,ik} \left(a_{\alpha \beta  }g^{\alpha \beta } \right)
  -\left(\lambda_{, \alpha \beta  } g^{\alpha \beta }\right) A_{ik}-\lambda_{, ki}   \left(A_{\alpha \beta} g^{\alpha \beta }\right)\right)  g_{jl}\\
&+&  \frac{1}{n}\left(-\left(\tau_{\alpha \beta } g^{\alpha \beta }\right)   g_{jl} +
 \left(g^{\alpha \beta }\Lambda_{,\alpha \beta}\right) a_{jl} +
 \Lambda_{,jl } \left(a_{\alpha \beta }g^{\alpha \beta } \right)
  -\left(\lambda_{, \alpha \beta } g^{\alpha \beta }\right) A_{jl}-\lambda_{, lj}   \left(A_{\alpha \beta } g^{\alpha \beta }\right)\right) g_{ik} \\
&=&
 \Lambda_{,jl} a_{ik}  +  \Lambda_{,ik} a_{jl}
  -\lambda_{, ki}   A_{jl}      -\lambda_{, lj} A_{ik}    . \end{array} $$

  Denoting $\left(\tau_{\beta\alpha } g^{\beta\alpha }\right)$ by $\tau$,
  $\left(\lambda_{,\beta\alpha } g^{\beta\alpha }\right)$ by $\mu$,
   $\left(\Lambda_{,\beta\alpha } g^{\beta\alpha }\right)$ by $\mathcal{M}$,   and using that
   $\left(a_{\alpha \beta  }g^{\alpha \beta } \right)=2 \lambda $, and  $\left(A_{\alpha \beta  }g^{\alpha \beta } \right)=2 \Lambda $,   one obtains

  $$\begin{array}{ll}
&\frac{1}{n}\left(-\tau   g_{ik} +
 \mathcal{M} a_{ik} +
 2\Lambda_{,ik}\lambda
  -\mu A_{ik}-2 \lambda_{, ki}  \Lambda \right)  g_{jl}\\
+& \frac{1}{n}\left(-\tau   g_{jl} +
 \mathcal{M}  a_{jl} +
2 \Lambda_{,jl } \lambda
  -\mu A_{jl}-2 \lambda_{, lj}   \Lambda \right) g_{ik} \\
=&
 \Lambda_{,jl} a_{ik}  +  \Lambda_{,ik} a_{jl}
  -\lambda_{, ki}   A_{jl}      -\lambda_{, lj} A_{ik}    . \end{array}
 $$

  Combining the terms, we obtain
\weg{   $$\begin{array}{ll}
&\frac{1}{n}\left(-\tau   g_{ik} +
 2 \Lambda_{,ik}\lambda
 -2 \lambda_{, ki}  \Lambda  \right)  g_{jl}
+ \frac{1}{n}\left(-\tau   g_{jl} +
 2 \Lambda_{,jl } \lambda
 -2 \lambda_{, lj}   \Lambda \right) g_{ik} \\
=& \left(-\frac{1}{n} \mathcal{M}    g_{jl} +
 \Lambda_{,jl}\right) a_{ik}  +
\left( \Lambda_{,ik}- \frac{1}{n}\mathcal{M}  g_{ik}  \right)   a_{jl}
  -\left(\lambda_{, ki}     -\frac{1}{n}\mu g_{ik}\right) A_{jl}    -\left( \lambda_{, lj}    -\frac{1}{n} \mu g_{jl}\right)A_{ik}  . \end{array} $$}

   \begin{equation}\begin{array}{ll}
&\left(-\frac{\tau}{2}   g_{ik} +
 2 \Lambda_{,ik}\lambda  \right)  \frac{g_{jl}}{n}  -\left(\frac{\tau}{2}   g_{ik}
 +\lambda_{, ki}  \Lambda  \right)  \frac{g_{jl}}{n}
+ \left(-\frac{\tau}{2}   g_{jl}
 -2 \lambda_{, lj}   \Lambda \right)\frac{g_{ik}}{n}  +
 \left(-\frac{\tau}{2}   g_{jl} +2  \Lambda_{,jl } \lambda  \right)\frac{g_{ik}}{n}   \\
=& \left(-\frac{1}{n} \mathcal{M}    g_{jl} +
 \Lambda_{,jl}\right) a_{ik}  +
\left( \Lambda_{,ik}- \frac{1}{n}\mathcal{M}  g_{ik}  \right)   a_{jl}
  -\left(\lambda_{, ki}     -\frac{1}{n}\mu g_{ik}\right) A_{jl}    -\left( \lambda_{, lj}    -\frac{1}{n} \mu g_{jl}\right)A_{ik}  . \end{array} \label{int18}\end{equation}

Now,  let us  calculate $\tau$: we multiply \eqref{int12} by $g^{il}$ and sum over  $i,l$. After dividing by $2$,   we obtain
  $$
 \tau
=\frac{2}{n}\left(
  \mathcal{M} \lambda
 - \mu   \Lambda \right) .$$

   Substituting in \eqref{int18}, we obtain

  $$\begin{array}{ll}
&\left(-\frac{\mathcal{M}}{n}     g_{ik} +
 \Lambda_{,ik} \right)  \frac{2\lambda}{n}g_{jl}  -\left(            -         \frac{\mu}{n}    g_{ik}
 +\lambda_{, ki}   \right)  \frac{2\Lambda}{n}g_{jl}
+ \left(\frac{\mu}{n}    g_{jl}
 -\lambda_{, lj}   \right)\frac{2\Lambda }{n}g_{ik}  +
 \left(-\frac{\mathcal{M}}{n}   g_{jl} + \Lambda_{,jl } \right)\frac{2\lambda }{n} g_{ik}  \\
=& \left(-\frac{1}{n} \mathcal{M}    g_{jl} +
 \Lambda_{,jl}\right) a_{ik}  +
\left( \Lambda_{,ik}- \frac{1}{n}\mathcal{M}  g_{ik}  \right)   a_{jl}
  -\left(\lambda_{, ki}     -\frac{1}{n}\mu g_{ik}\right) A_{jl}    -\left( \lambda_{, lj}    -\frac{1}{n} \mu g_{jl}\right)A_{ik}  , \end{array} $$
 which is equivalent to

  \begin{equation}\begin{array}{cl}
&\left( \Lambda_{,ik}-\frac{\mathcal{M}}{n}     g_{ik}  \right)
\left(a_{jl} -\frac{2\lambda }{n}g_{jl}\right) +\left( \Lambda_{,jl }-\frac{\mathcal{M}}{n}   g_{jl}  \right)\left(a_{ik}-\frac{2 \lambda  }{n} g_{ik} \right)\\
 = & \left(           \lambda_{, ki}  -         \frac{\mu}{n}    g_{ik}    \right)
 \left(  A_{jl} -\frac{2\Lambda}{n}g_{jl}\right)+
 \left(\lambda_{, lj} -\frac{\mu}{n}    g_{jl}   \right) \left(A_{ik} -\frac{2 \Lambda  }{n}g_{ik} \right)
 \end{array} \label{int20}\end{equation}

 Denoting $B_{ik}:= \left(\Lambda_{,ik}-\frac{\mathcal{M}}{n}     g_{ik}
  \right) $ ,     $ b_{ik}:=  \left(\lambda_{, ki}
 - \frac{\mu}{n}    g_{ik} \right) $,
 $d_{jl}:=   \left(a_{jl} -\frac{2\lambda }{n}g_{jl}\right) $,
 $D_{jl}:= \left(A_{ik} - \frac{2\Lambda  }{n}g_{ik}\right)$, we see that \eqref{int20} is equivalent to
 \begin{equation}
B_{ik} d_{jl} + B_{jl} d_{ik} = b_{ik} D_{jl} + b_{jl} D_{ik}.
 \label{int21}\end{equation}

Now,  it is easy to see that the condition $B_{\alpha} d_{\beta} + B_{\beta} d_{\alpha} = b_{\alpha} D_{\beta} + b_{\beta} D_{\alpha}$ for dimensions $\ge 3$ implies
that $B$ is proportional to $b$ and $D$ to $d$, or that  $B$ is proportional to $D$ and $b$ to $d$.
We see that the condition \eqref{int21} is essentially the same as this  condition, the role of $\alpha $
and  $\beta$ play the multi-indexes $ik$ and $jl$. Since $g,$ $  a,$ and $A$ are linearly independent,
 $D$ can not be  proportional to $d$. Thus,   $b$ is proportional to $d$ which imply that
 $\lambda_{,ij}$
 is a linear combination of  $a_{ij}$  and $g_{ij} $. The coefficients of the linear combination are not important for us, let
 $$ \lambda_{,ij} = \rho \cdot g_{ij} +  K\cdot a_{ij}$$

 Substituting this equation in \eqref{int1}, we obtain

\begin{equation}  a_{i \alpha} Z^{\alpha}_{jkl} +  a_{\alpha j}Z^{\alpha}_{ikl}=0,  \label{int32}
\end{equation}
 Where $Z^{i}_{jkl}= R^i_{jkl} - K\cdot  ( \delta^i_lg_{jk} - \delta^i_kg_{jl}  )$.
Clearly, the equation \eqref{int32} is equivalent to  the equation \eqref{int33b}.

We see that by construction the tensor $Z_{ijkl}$ is skew-symmetric with respect to the first two indexes
(actually, the tensor $Z$ has the same symmetries as the curvature tensor, i.e., for example,  skew-symmetric with respect to the last  two indexes   as well,   but we will need only the first  two indexes).  Thus,
Lemma~\ref{13} and, therefore,  Theorem~\ref{thm1}  are proved, \qed

\section{Fubini theorem and  sectional operators  \\ on semisimple Lie algebras}

In this section we discuss an unexpected and remarkable relationship between geodesically equivalent  metrics and some special operators on semisimple Lie algebras which appeared in the theory of integrable systems.

We start with a brief  overview on (one special type of)  integrable Euler equations on semisimple Lie algebras (see \cite{bogo, fomtrofbook, marsden, MF1, pere,  reyman}  for details).

Let  $\goth g$ be a semisimple Lie algebra, $R: \goth g \to \goth g$  an operator  symmetric with respect to the Killing form $\langle~,~\rangle$ on $\goth g$. The differential equation
\begin{equation}
\dot x = [R(x),x], \quad x\in\goth g,
\label{Euler}
\end{equation}
is Hamiltonian on $\goth g$ with respect to the standard Lie-Poisson structure  and called the {\it Euler equation} related to the Hamiltonian function $H(x)=\frac{1}{2}\langle R(x), x\rangle$.

A classical, interesting  and extremely difficult problem is to find those operators
$R: \goth g \to \goth g$ for which the system (\ref{Euler}) is completely integrable.

One of such operators was discovered by S.\,Manakov in \cite{Manakov} and his idea then led to an elegant general construction developed by A.\,Mischenko and A.\,Fomenko \cite{MF1} and  called {\it argument shift method}. This construction in brief can be presented as follows.

Assume that $R : \goth g \to \goth g$ satisfies the following identity
\begin{equation}
[R(x),a]=[x,b],   \quad x\in \goth g,
\label{sectional}
\end{equation}
for certain $a,b\in \goth g$,  $a\ne 0$.   Then the following statement holds

\begin{Th}{\rm  \cite{MF1}} \label{MFth}
Let $R: \goth g \to \goth g$ be symmetric and satisfy {\rm (\ref{sectional})}. Then

1) the system {\rm (\ref{Euler})}
admits the following Lax representation with a parameter:
$$
 \frac{d}{dt}(x + \lambda a) = [ R(x)+\lambda b, x+\lambda a];
$$

2) the functions  $f(x+\lambda a)$, where $f:\goth g\to \R$ is an invariant of the adjoint representation,  are first integrals of {\rm (\ref{Euler})}  for any $\lambda\in \R$ and, moreover, these integrals commute;

3) if $a\in \goth g$ is  regular, then {\rm (\ref{Euler})} is completely integrable.
\end{Th}

This construction has a very important particular case. If the Lie algebra $\goth g$ admits  a  $\Z_2$-grading,  i.e., a decomposition $\goth g = \goth h  + \goth v$  (direct sum of  subspaces) such that $[\goth h,\goth h]\subset \goth h$, $[\goth h,\goth v]\subset \goth v$, $[\goth v,\goth v]\subset \goth h$,
then  we may consider  $R:\goth h \to \goth h$  satisfying   (\ref{sectional}) with $a,b\in \goth v$,  and
Theorem   \ref{MFth} still holds if we replace $\goth g$ by $\goth h$.

The most important  example for applications (in particular,  in the theory of integrable tops)  is $\goth  g = sl(n,\R)$,  $\goth h = so(n,\R)$, with  $a$ and $b$ symmetric matrices. This  is  the  situation  that was studied in the pioneering work by S.~Manakov \cite{Manakov}  leading  to integrability of the Euler equations of $n$-dimensional rigid body dynamics.

From the algebraic point of view,  the above construction still makes sense if we replace $so(n)$ by $so(p,q)$ and  assume $a,b$ to be  symmetric operators with respect to the corresponding indefinite form $g$.
Moreover, if  we complexify our considerations we do not even notice any difference.
However, to indicate the presence (but not influence) of the bilinear form $g$, we shall denote the space of  $g$-symmetric operators by $Sym(g)$, and the Lie algebra of $g$-skew-symmetric operators by $so(g)$.

\begin{Def}
\label{defsect} 
{\rm  We shall say that $R: so(g) \to so(g)$ is a {\it sectional operator} associated with $a,b\in Sym(g)$,  if  $R$ is symmetric with respect to the Killing form and the following identity holds:
\begin{equation}
[R(x), a] = [x,b],    \quad \mbox{for all } x\in  so(g).
\label{sectionall}
\end{equation}
}\end{Def}

We follow the terminology introduced by Fomenko and Trofimov  in \cite{troffom, fomtrofbook} where they studied various generalizations of such operators. Strictly speaking, the above definition is just a particular case of  a more general construction.  The term``sectional"  was motivated by the following reason.   The  identities
(\ref{sectional}) and (\ref{sectionall}) suggest  that one may represent $R$ as $\mbox{ad}_a^{-1}\mbox{ad}_b$,  but in general we cannot do so because $\mbox{ad}_a$, as a rule,  is not invertible. That is why  the operator $R$ splits into parts each of which acts independently on its own subspace (section).

A surprising relationship between sectional operators and geodesically equivalent metrics is explained by the following observation. Notice, first of all, that due to its algebraic symmetries   (skew-symmetry with respect to $i,j$ and $k,l$  and symmetry with respect to permutation of pairs $(ij)$ and $(kl)$),   the Riemann curvature tensor $R_{ij,kl}$  can be naturally considered as
a symmetric operator  $R: so(g) \to so(g)$   (strictly speaking we need to raise indices $i$ and $k$  by means of $g$ to get the tensor of the form $R^{i,k}_{j,l}$). Having this interpretation of $R$ in mind, we immediately obtain

\begin{Th}
Let  $g$ and $\bar g$ be geodesically equivalent  nonproportional metrics, then
the Riemann curvature tensor $R^{i,k}_{j,l}$ of the metric $g$ is a sectional operator in the sense of Definition {\rm \ref{defsect}}.
More precisely,
$$
[R(x), a]=[x,b]
$$
where $a$ is  the $g$-symmetric operator associated with the form $a_{ij}$ defined by (3),
and  $b$ is the $g$-symmetric operator associated with the form $2\lambda_{,ij}$ (Hessian  of  $\mathrm{tr}\, a$).
\end{Th}

The proof of this statement is just the observation that (\ref{sectionall}) is a translation of the compatibility condition (\ref{int1}) into Lie-algebraic language.

Before discussing the proof of  the Fubini theorem in this "new" language,  we make some remarks which could also be useful.

Notice, first of all, that in our new notation the condition "curvature is constant"   (at a point) simply means that  $R: so(g) \to so(g)$ is a scalar operator, i.e. $R(x)=K\cdot x$.

Furthermore, it is a very simple Lie-algebraic fact that   (\ref{sectionall})   implies that $a$ and $b$ commute.  Indeed,
$\langle [b,a], x \rangle = \langle a, [x,b] \rangle =\langle a, [R(x),a] \rangle =
\langle [a,a], R(x) \rangle =0 $ for any $x\in so(g)$,  so $[a,b]=0$.   In the theory of projectively equivalent metrics this means that  the operator $a$ commutes with the Hessian of its trace $2\lambda =\tr\, a$. This fact is, of course,  well known, but  the above proof seems to be the simplest one.  Moreover, if instead of $a$ we substitute any element  $\xi$  from its centralizer $C(a)$,  we obviously get the same conclusion
$[b,\xi]=0$, i.e., $b$  lies in the center of the centralizer of $a$.  This means, in fact, that  $b$ is a polynomial of $a$.

Finally, if $a$ is regular in the Lie-algebraic sense, i.e. its minimal polynomial coincides with the characteristic one, then  the operator ${\rm ad}_a: so(g) \to Sym(g)$  has trivial kernel  so that  the sectional operator $R$  (i.e., the curvature tensor!) can be reconstructed from $a$ and $b$.   Namely, $R(x)={\rm ad}_a^{-1}{\rm ad}_b (x)$, a well-known formula in the theory of integrable systems on Lie algebras.  If we take into account the fact that  $b=P(a)= \lambda_{n-1} a^{n-1} + \lambda_{n-2} a^{n-2} + \dots + \lambda_1 a + \lambda_0$  (polynomial of $a$), then this formula can be rewritten as
\begin{equation}
R(x) = \frac{d}{dt}\left. P(a+tx) \right|_{t=0} .
\label{rx}
\end{equation}
Indeed,  $[P(a+tx), a+tx]=0$ implies
$$
0=\frac{d}{dt} [P(a+tx), a+tx]|_{t=0}=[\frac{d}{dt} P(a+tx) |_{t=0},a]+[P(a), x],
$$
i.e., $[\frac{d}{dt} P(a+tx)|_{t=0} ,a]=[x,b]$. Since $a$ is regular, we have (\ref{rx}).

This shows, in particular,  that the algebraic structure of the curvature tensor can be understood in terms of the operator $a$ only.

Using this language we now give another proof of  the tensor part of  the Fubini theorem  (Lemmas 2 and 3).

Assume that we have three geodesically equivalent metrics  $g$, $\bar g$, and $\hat g$.
Then the Riemann curvature tensor $R$ of the metric $g$ satisfies at the same time two identities :
\begin{equation}
[ R(x), a]=[x, b]\qquad \mbox{and}  \qquad [ R(x), A]=[x, B],
\label{2sectional}
\end{equation}
where $a^k_j = g^{ki}a_{ij}$, $A^k_j = g^{ki}A_{ij}$, $b^k_j = 2g^{ki}\lambda_{, ij}$, $B^k_j = 2g^{ki}\Lambda_{,ij}$ (cf. (\ref{int1}) and (\ref{int2})).

From now on,  we may forget about the geometrical meaning of $a,b, A, B$ and start thinking of them as just certain $g$-symmetric operators.  In addition, without loss of generality we may assume all these operators to be trace free (as, of course, it should be in the semisimple Lie algebra $sl(n,\R)$ which stands behind this construction). Moreover,  we are allowed to complexify all the objects  so that instead of  $so(g)$ and $Sym(g)$  we may simply consider the spaces of symmetric and skew-symmetric complex matrices.

The reformulation of the (algebraic part of) Fubini theorem are the following  analogs of Lemmas 2 and 3 respectively.

\begin{Lemma}
Let $R: so(g) \to so(g)$  be symmetric and satisfy {\rm (\ref{2sectional})}.
If $a$ and $A$  are not proportional,  then
$b$ is proportional to $a$  and, therefore,  $[R(x) - K\cdot x, a]=0$ for some $K\in \R$.
\label{lemma4}
\end{Lemma}

\begin{Lemma}
If  $a$ is regular, i.e., its minimal polynomial coincides with the characteristic one, then
the identity $[R(x) - K\cdot x, a]=0$   implies  $R = K\cdot {\rm{id}}$  (i.e., the curvature is constant).
\label{lemma5}
\end{Lemma}

{\bf Proof of Lemma \ref{lemma4}.}
Let $y$ and $z$ be arbitrary $g$-symmetric matrices, then $[A, y], [a, z] \in so(g)$
and we have:
$$
[ R([A, y]), a]=[[A,y], b], \quad  [ R([ a , z]), A]=[[a, z], B].
$$
Since $R$ is symmetric with respect to the Killing form $\langle \ , \ \rangle$  we
have
$$
\langle [[A, y], b] , z\rangle= \langle [ R([A, y]), a], z\rangle =
\langle  R([A, y]), [ a , z] \rangle = \langle  [A, y], R([ a , z]) \rangle =
$$
$$
 \langle  y, [ R([ a , z]), A] \rangle =
\langle  y , [[a, z], B] \rangle = \langle [[B,y], a], z \rangle
$$
Since $z$ is an arbitrary symmetric matrix, we conclude that
\begin{equation}
[[A, y], b] =[[B,y], a].
\label{mainident}
\end{equation}

This relation is an analog of (\ref{int8}). 
Similarly,  $ [[a, y], B] =[[b,y], A]$.
Using the Jacobi identity, it is not hard to see that
$$
[b, A]=[a, B]
$$

Rewriting (\ref{mainident}) as
$$
y (Ba - A b) + (aB - b A) y = Bya+ ayB - byA - Ayb
$$
and noticing that $[b, A]=[a, B] $ implies $Ba - Ab
=aB - b A$,  we get
$$
y  T + T  y = Bya+ ayB - byA - Ayb
$$
where $T$ denotes $aB - b A$  (this is an analog of $\tau$ from (\ref{int10})).

This formula can be considered as a relation between two  linear operators acting on the space of symmetric matrices  (the argument of both operators is $y\in Sym(g)$).
To get some consequences from this identity,  we take a kind of its trace. Recall that we consider $A, a, B, b, y, T$ as usual  symmetric (complex) matrices.

Instead of $y$ we substitute the symmetric matrix of the form  $e_i  v^\top + v e_i^\top$,  where $e_i$ and $v$ are vector-columns ($e_1,\dots, e_n$  is the standard (orthonormal) basis), then apply the result to  $e_i$ and take the sum over $i$. Here is the result:
$$
(e_i  v^\top + v e_i^\top) T e_i + T ( e_i  v^\top + v e_i^\top) e_i=
B ( e_i  v^\top + v e_i^\top)a e_i + ...
$$
$$
e_i  (Tv, e_i) +  v (Te_i, e_i)+ Te_i (v,e_i) + Tv (e_i, e_i) =
B e_i (a v, e_i) + B v (a e_i, e_i) + ...
$$

Using obvious facts from Linear Algebra such as
$$
\sum_i(Te_i, e_i) = \tr \,T,  \quad  \sum_i(e_i,e_i)=n, \quad
\sum_i e_i (v, e_i)=v,
$$
we get
$$
Tv +  \tr \, T\cdot v  + Tv + n \cdot Tv  =
B a v  +  \tr \, a\cdot B v  + ...
$$
Taking into account that $a, A, b, B$ are all trace free we have
$$
((n+2) T +   \tr\, T \cdot {\rm Id} ) v  =
(B a   + aB -   A b - bA) v
$$
Since $v$ is arbitrary and $T=B a -   A b = aB  - bA$, we finally get
$$
nT + \tr \, T \cdot {\rm Id }= 0
$$

But this simply means that $T=0$. Hence we come to the identity of the form
\begin{equation}
Bya+ ayB = byA + Ayb.
\label{mainident2}
\end{equation}

It remains to use the following simple statement:
if  $a,  b,  A, B$  are symmetric, $a\ne 0$ and (\ref{mainident2}) holds for any symmetric $y$,
then either    $b = K\cdot  a$, or  $A = K\cdot a$ for some constant $K\in \R$. 

By our assumption, $a$ and $A$ are not proportional, so we conclude that
$b= K\cdot a$ and therefore the identity $[ R(x), a]=[x, b]$ becomes $[ R(x)- K\cdot x, a]=0$, as needed.

Notice that  (\ref{mainident2}) and  the rest of the proof almost literally repeat (\ref{int21}) 
and the end of the proof of Lemma 2, \qed

{\bf Proof of Lemma \ref{lemma5}.} 
Let $[ R(x)- K\cdot x, a]=0$ and $a$ be regular.  It is a well known algebraic fact that the centralizer of a regular matrix $a$ is generated by the powers of $a$. In particular, the centralizer of $a$ consists of $g$-symmetric matrices.  On the other hand, $R(x)- K\cdot x$ is skew-symmetric. Thus, $R(x)- K\cdot x$ has to be zero for any $x$, i.e., $R = K\cdot {\rm{id}}$, as was to be proved, \qed

\label{bolsinov}

\vskip 8pt

Alexey V.Bolsinov, \newline
School of Mathematics, Loughborough University, \newline
Loughborough, LE11 3TU, UK  \newline
e-mail: \quad A.Bolsinov@lboro.ac.uk,

Volodymyr  Kiosak, \newline
Institute of Mathematics, FSU Jena, \newline 
07737 Jena Germany  \newline
e-mail: \quad kiosak@minet.uni-jena.de

Vladimir S. Matveev\newline
Institute of Mathematics, FSU Jena, \newline 
07737 Jena Germany  \newline
e-mail: \quad matveev@minet.uni-jena.de


\begin{thebibliography}{99}
\bibitem{Aminova1}
A. V. Aminova,
\emph{A Lie problem, projective groups of two-dimensional
Riemann surfaces, and solitons},
Izv. Vyssh. Uchebn. Zaved. Mat. 1990, no. 6, 3--10;
translation in Soviet Math. (Iz. VUZ) \textbf{34} (1990), no. 6, 1--9.

\bibitem{Aminova} A. V.  Aminova, {\it Pseudo-Riemannian manifolds with general geodesics,}   Russian Math. Surveys  48  (1993),  no. 2, 105--160, MR1239862, Zbl 0933.53002.


\bibitem{Aminova2}  A. V. Aminova,
\emph{Projective transformations of pseudo-Riemannian manifolds. Geometry, 9.}
 J. Math. Sci. (N. Y.) \textbf{113} (2003), no. 3, 367--470.
\bibitem{Benenti3} S. Benenti, \emph{ Special symmetric two-tensors, equivalent dynamical systems, cofactor and bi-cofactor systems,}
Acta Appl. Math. {\bf 87}(2005), no. 1-3, 33--91.

 \bibitem{Beltrami} E. Beltrami,
{\it
Resoluzione del problema: riportari
i punti di una
superficie sopra un piano in modo che le linee geodetische vengano
rappresentante da linee rette},
Ann. Mat., {\bf 1}(1865), no. 7, 185--204.


\bibitem{bogo}
O. I. Bogoyavlensky, \emph{
Integrable Euler equations on Lie algebras, arising in problems of mathematical physics.} (Russian) 
Izv. Akad. Nauk SSSR Ser. Mat. \textbf{48} (1984), no. 5, 883--938. MR0764304 



\bibitem{benenti} A.V.~Bolsinov, V.S.~Matveev,
{\it Geometrical interpretation of Benenti systems},
J. of Geometry and Physics,  {\bf 44}(2003),  489--506, MR1943174, Zbl 1010.37035.

\bibitem{bryant}   R. L. Bryant,
G. Manno,
V.  S. Matveev, \emph{A solution of a problem of Sophus Lie: Normal forms of 2-dim metrics admitting two projective vector fields}, accepted to Math. Ann,
    arXiv:0705.3592 .


\bibitem{duna}   R. L. Bryant, M. Dunajski, M. Eastwood,
\emph{ Metrisability of two-dimensional projective structures, }    arXiv:0801.0300.



\bibitem{eastwood} M. Eastwood,  V. S. Matveev,  \emph{ Metric connections in projective differential geometry,}
 Symmetries and Overdetermined Systems of Partial Differential Equations (Minneapolis, MN, 2006), 339--351,
 IMA Vol. Math. Appl.,    {\bf
144}(2007),   Springer, New York.

\bibitem{fomtrofbook}
A.T.Fomenko, V.V.Trofimov, \emph{Integrable systems on Lie Algebras and Symmetric Spaces},  Gordon and Breach, 1988.

\bibitem{Fubini1} G. Fubini, {\it Sui gruppi transformazioni geodetiche, }
Mem. Acc. Torino  {\bf 53}(1903), 261--313.
\bibitem{Fubini2} G. Fubini, {\it Sulle coppie di varieta geodeticamente applicabili,}  Acc. Lincei {\bf 14}(1905),
 678--683 (1°Sem.),
315--322 (2\ °Sem.).

\bibitem{kiosak} V. A. Kiosak, \u I. Mikesh, \emph{ On the degree of mobility of Riemannian spaces with respect to geodesic mappings. }  (Russian)  The geometry of imbedded manifolds (Russian),  35--39, {\bf 124}, Moskov. Gos. Ped. Inst., Moscow, 1986.

 \bibitem{kiosak-mikes} V. A. Kiosak, \u I. Mikesh, \emph{ On geodesic mappings of Einstein spaces.} (Russian)  Izv. Vyssh. Uchebn. Zaved. Mat.  2003, , no. 11, 36--41;  translation in  Russian Math. (Iz. VUZ)  {\bf 47 } (2003),  no. 11, 32--37 (2004)
 
\bibitem{koenigs}
G. Koenigs,
\emph{Sur les g\'eodesiques a int\'egrales quadratiques, }
Note II from Darboux'
`Le\c{c}ons sur la th\'eorie g\'en\'erale des surfaces',
 Vol. IV, Chelsea Publishing, 1896.

 \bibitem{era1} B. S. Kruglikov,  V. S. Matveev, {\it Vanishing of the entropy pseudonorm for certain integrable systems, }  Electron. Res. Announc. Amer. Math. Soc.  {\bf 12}  (2006), 19--28.

 \bibitem{kruglikov} B. Kruglikov, \emph{
Invariant characterization of Liouville metrics and polynomial integrals,}  arXiv:0709.0423
\bibitem{lagrange} J.-L.  Lagrange, \emph{ Sur la construction des cartes g\'eographiques,} Nov\'eaux M\'emoires de l'Acad\'emie des Sciences et Bell-Lettres de Berlin, 1779.

 \bibitem{Levi-Civita}
 T. Levi-Civita, {\it Sulle trasformazioni delle equazioni
 dinamiche}, Ann. di Mat., serie $2^a$, {\bf 24}(1896), 255--300.

\bibitem{Lie}
S. Lie,
\emph{Untersuchungen \"uber geod\"atische Kurven},
Math. Ann. \textbf{20} (1882); Sophus Lie Gesammelte Abhandlungen,
Band 2, erster Teil, 267--374. Teubner, Leipzig, 1935.

\bibitem{Manakov}
 S.V. Manakov, \emph{Note on the integration of Euler's equation of the dynamics of an $N$- dimensional rigid body}, Funct. Anal. Appl. \textbf{11} (1976), 328-329.

\bibitem{marsden}
J.E.~Marsden, T.S.~Ratiu:
\emph{Introduction to Mechanics 
and Symmetry.}
Springer Verlag, New York, 1999. 


\bibitem{MT}
V. S. Matveev, P. J. Topalov,
\emph{Trajectory equivalence and corresponding integrals},
 Regular and Chaotic Dynamics,
\textbf{3} (1998), no. 2, 30--45.

\bibitem{dim2}
V. S. Matveev, P. J. Topalov,
\emph{Geodesic equivalence of metrics on surfaces, and their
integrability,} Dokl. Math. \textbf{60} (1999), no.1, 112-114.

{\bibitem{ERA}
V. S. Matveev and P. J. Topalov,
\emph{Metric with ergodic geodesic flow is completely determined by
unparameterized geodesics,} ERA-AMS, \textbf{6} (2000), 98--104.}

 \bibitem{quantum}  V. S. Matveev, P. J. Topalov,
 {\it Quantum integrability for the   Beltrami-Laplace
 operator as geodesic  equivalence,}
  Math. Z. {\bf 238}(2001), 833--866, MR1872577, Zbl 1047.58004.


\bibitem{integrable} V. S. Matveev and  P. J. Topalov,
 {\it Integrability  in theory of geodesically equivalent metrics},
J. Phys. A., {\bf 34}(2001), 2415--2433, MR1831306, Zbl 0983.53024.



\bibitem{topology} V. S. Matveev,
  {\it Three-dimensional manifolds having metrics with the same geodesics,}
    Topology     {\bf 42}(2003) no. 6,  1371-1395, MR1981360, Zbl 1035.53117.


\bibitem{hyperbolic}  V. S. Matveev,
{\it Hyperbolic manifolds are geodesically rigid,}
 Invent. math. {\bf 151}(2003), 579--609, MR1961339, Zbl 1039.53046.

\bibitem{obata}
V. S. Matveev,
\emph{Die Vermutung von Obata f\"ur Dimension $2$,}
Arch. Math. \textbf{82} (2004), 273--281.

\bibitem{sol}
V. S. Matveev,
\emph{Solodovnikov's theorem in dimension two},
Dokl. Math. \textbf{69} (2004), no. 3, 338--341.

\bibitem{CMH}
V. S. Matveev,
\emph{Lichnerowicz-Obata conjecture in dimension two,}
Comm. Math. Helv. \textbf{81}(2005) no. 3,  541--570.

\bibitem{fomenko}  V. S. Matveev,   { \it  Beltrami problem, Lichnerowicz-Obata conjecture
and  applications of integrable systems in differential geometry,} Tr. Semin.
  Vektorn.   Tenzorn. Anal, {\bf 26}(2005), 214--238.

\bibitem{beltrami_short}
V. S. Matveev,
\emph{Geometric explanation of Beltrami theorem},
Int. J. Geom. Methods Mod. Phys. \textbf{3} (2006), no. 3, 623--629.

\bibitem{sbornik}  V. S. Matveev, \emph{ On degree of mobility of complete metrics,}
Adv. Stud. Pure Math., {\bf 43}(2006), 221--250.

\bibitem{archive}
V. S. Matveev,
\emph{Proof of projective Lichnerowicz-Obata conjecture},
 J. Diff. Geom. (2007), {\bf 75}(2007),  459--502.

\bibitem{alone} V. S. Matveev, {\it A solution of another  S. Lie Problem: 2-dim metrics admitting  projective vector field,}   Math. Ann.,  submitted. arXiv:math/0802.2344


\bibitem{conformal} V. S. Matveev, H.-B.  Rademacher, M. Troyanov, and A. Zeghib,  {\it Conformal Lichnerowicz-Obata conjecture, }   Math. Res. Let.,  submitted.
arXiv:math/0802.3309

\bibitem{mikes}
 J. Mikes, {\it
  Geodesic mappings of affine-connected and Riemannian spaces.
Geometry, 2.,}
  J. Math. Sci. {\bf 78}(1996), no. 3, 311--333.


\bibitem{MF1}
A. S. Mishchenko, A. T. Fomenko,
{Euler equations on finite-dimensional Lie groups}.
Izv. Acad. Nauk SSSR, Ser. matem. {\bf 42}, no.2,  396-415  (1978)
(Russian);    English translation: Math. USSR-Izv. {\bf 12}, no.2,
371-389  (1978).




\bibitem{pere}
A. M. Perelomov, \emph{ Integrable systems of classical mechanics and Lie algebras. Vol. I.} BirkhŠuser Verlag, Basel, 1990, MR1048350. 


\bibitem{reyman}
A.G.~Reyman,  M.A.~Semenov-Tian-Shansky,  \emph{Integrable systems. A group theoretic approach.}
Institute of Computer Science, Moscow-Izhevsk, 2003 (Russian).



\bibitem{schouten}
J. A. Schouten,
\emph{Erlanger Programm und \"Ubertragungslehre.
 Neue Gesichtspunkte zur Grundlegung der Geometrie},
Rendiconti Palermo  {\bf 50} (1926), 142--169.



\bibitem{schur} F. Schur, {\it Ueber den Zusammenhang der R\"aume constanter
Riemann'schen Kr\"ummumgsmaasses mit den projektiven R\"aumen, } Math. Ann. {\bf 27}(1886), 537--567.

{\bibitem{shandra} I. G. Shandra, \emph{  On the geodesic mobility of Riemannian spaces,}  Math. Notes {\bf 68}(2000), no. 3-4, 528--532.}


\bibitem{Sinjukov} N. S.  Sinjukov,
 {\it Geodesic mappings of Riemannian spaces},  (in Russian)
``Nauka'', Moscow, 1979, MR0552022, Zbl 0637.53020.


\bibitem{s1} A. S. Solodovnikov, {\it Projective transformations of Riemannian spaces,}
 Uspehi Mat. Nauk (N.S.) {\bf 11}(1956), no. 4(70), 45--116, MR0084826, Zbl 0071.15202.

\bibitem{s2} A. S.  Solodovnikov, {\it
 Spaces with common geodesics,} Trudy Sem.
 Vektor. Tenzor. Anal. {\bf 11}(1961), 43--102, MR0163257, Zbl 0161.18904.

\bibitem{s3} A. S. Solodovnikov, {\it  Geometric description of all possible representations of a Riemannian metric in Levi-Civita form,}
 Trudy Sem. Vektor. Tenzor.
  Anal. {\bf 12}(1963), 131--173, MR0162201, Zbl 0163.43302.


\bibitem{involutivity} P. Topalov, \emph{
 Geodesic hierarchies and involutivity,}
J. Math. Phys.  {\bf 42}(2001), no. 8, 3898--3914.


\bibitem{Topalov} P. Topalov, {\it  Commutative conservation laws for geodesic flows of metric admitting projective symmetry,}
 Math. Res. Lett. {\bf 9}(2002), no. 1, 65--72, MR1892314, Zbl 1005.53065.

\bibitem{dedicata}
P. J. Topalov and V. S. Matveev,
\emph{Geodesic equivalence via integrability},
Geometriae Dedicata \textbf{96} (2003), 91--115.


\bibitem{troffom}
V. V. Trofimov and A. T. Fomenko,
\emph{Dynamical systems on the orbits of linear representations of Lie groups and the complete integrability of certain hydrodynamical systems},
Functional Analysis and Its Applications  \textbf{17} (1983), 2--29.


\bibitem{Vries} H. L. de  Vries, {\it \"Uber Riemannsche R\"aume, die
infinitesimale konforme Transformationen gestatten,}
 Math. Z. {\bf 60}(1954), 328--347, MR0063725, Zbl 0056.15203.


\bibitem{Weyl}  H. Weyl, {\it Zur Infinitisimalgeometrie: Einordnung der projektiven
und der  konformen Auffasung,} Nachrichten von der K. Gesellschaft
der Wissenschaften zu G\"ottingen, Mathematisch-Physikalische
Klasse, 1921;
 ``Selecta Hermann Weyl'', Birkh\"auser Verlag,
   Basel und Stuttgart,
1956.

\bibitem{zaharov} V.E.  Zakharov, {\it Description of the $n$-orthogonal curvilinear coordinate systems and Hamiltonian integrable systems of hydrodynamic type. I. Integration of the Lam\'e equations.}  Duke Math. J.  {\bf 94}(1998),  no. 1, 103--139.

\bibitem{zakirova} Z. Zakirova,  \emph{ On projective group properties of the 6D pseudo-Riemannian space,}   Czechoslovak J. Phys.  55  (2005),  no. 11, 1541--1544.

\end{thebibliography}
\end{document}